\documentstyle [11pt,amstex,psfig]{article}
\def\al{\alpha}

\def\text#1{{\em #1}}

\def\ep{\varepsilon}

\def\be{\begin{equation}}
\def\ee{\end{equation}}
\def\bear{\begin{eqnarray}}
\def\eear{\end{eqnarray}}
\def\best{\begin{eqnarray*}}
\def\eest{\end{eqnarray*}}
\def\pf{{\bf Proof}: }

\newtheorem{theorem}{Theorem}[section]
\newtheorem{prop}[theorem]{Proposition}
\newtheorem{proposition}[theorem]{Proposition}
\newtheorem{lemma}[theorem]{Lemma}
\newtheorem{cor}[theorem]{Corollary}

\newtheorem{defn}[theorem]{Definition}

\def\rem{ \addtocounter{theorem}{1}
{\non \bf Remark \arabic{section}.\arabic{theorem} }}

\def\non{\noindent}
\def\pf{\non {\bf Proof. }}
\def\qed{\nopagebreak \hskip .1in { $\Box$ }\penalty10000 %
\hskip\parfillskip \par  }

\def\ra{\rightarrow}

\def\r#1{\right#1}
\def\l#1{\left#1}

\def\ma#1{\mathop {#1} \limits}

\def\Si{\Sigma}
\def\De{\Delta}

\def\ti{\times}

\def\Z{{ \Bbb Z}}
\def\R{{ \Bbb R}}
\def\P{{ \Bbb P}}
\def\Q{{ \Bbb Q}}
\def\cx{{ \Bbb C}}

\def\ev{\mbox{\rm ev}}

\def\ov#1{\overline{#1}}

\def\w{\omega}
\def\tr{\mbox{tr}}
\def\M{{\cal M}}

\title{\bf Gromov Invariants and Symplectic Maps\vskip.2in}
\author{ Eleny-Nicoleta Ionel\thanks{partially supported by a M.S.R.I.
Postdoctoral Fellowship} \\  M.I.T.\\
 Cambridge, MA  02139 \and Thomas H. Parker\thanks{partially supported by
N.S.F. grant
DMS-9626245}\\ Michigan State University\\ East Lansing, MI   48824}
\date{March 1, 1997}

\addtocounter{section}{-1}
\begin{document}

\maketitle

\ \bigskip
\setcounter{equation}{0}
\section{\bf Introduction}
\bigskip

Given a symplectomorphism $f$ of a symplectic manifold $X$, one can form
the `symplectic mapping cylinder'
\begin{equation}
X_f = {X\times \R\times S^1}/\Z
\label{eq00.1}
\end{equation}
where the $\Z$ action is generated by $(x,s,\theta)\mapsto
(f(x),s+1,\theta)$.  In
this paper we compute the Gromov invariants of the manifolds $X_f$ and of
 fiber sums of the $X_f$ with other symplectic
manifolds.  This is done by expressing the Gromov invariants in terms of
the Lefschetz zeta function of
$f$ and, in special cases, in terms of the Alexander polynomials of knots.
The result is a large set of
interesting non-K\"{a}hler symplectic manifolds with computational ways of
distinguishing them.  In
particular, this gives a simple  symplectic construction of the
`exotic' elliptic surfaces recently discovered by Fintushel and Stern and
of related `exotic' symplectic
6-manifolds.

\medskip
 A closed symplectic manifold $(X,\w)$ has three `classical invariants':

\medskip

\hskip.5in  (i) the diffeomorphism type of $X$,
\smallskip

\hskip.5in  (ii) the DeRham cohomology class of $\w$, and
\smallskip

\hskip.5in  (iii) the homotopy class $[J]$ of the almost complex structures
$J$ compatible with $\w$.

\medskip

\noindent ($[J]$ determines the canonical class $\kappa\in H^2(X;\Z)$).
For many purposes these invariants are inadequate for distinguishing
symplectic manifolds.
In 1985 Gromov proposed constructing new symplectic invariants by
introducing  a $J$, considering the
$J$-holomorphic curves in $X$, and mimicking the construction of the
enumerative invariants of algebraic
geometry.  After much development, there now exist four or five related but
distinct versions of these ``Gromov
invariants'' ([MS], [RT], [LT],  [T3]).

In section 1 we define the ``degree zero Gromov invariants'' that we will
use.  These are built from the
invariants of Ruan and Tian, which count perturbed holomorphic curves
([RT]).  More specifically, for each class
$A\in H_2(X)$, we assemble Ruan-Tian invariants into a `partial Gromov
series' $Gr^A(X)$; this is a power series
in a variable $t_A$ whose coefficients  count, in a rather non-obvious way,
the  perturbed holomorphic curves
representing multiples of $A$.  The full degree zero Gromov series
$Gr(X)$ is the product of the $Gr^A(X)$ over
all primitive classes $A$.  These invariants are defined in all
dimensions.  In dimension four they can sometimes be related to the
Seiberg-Witten invariants by applying
the Theorem of Taubes [T1] and a result of the authors [IP1].

\bigskip

Now fix a  symplectomorphism $f$ of a closed symplectic manifold $X$, and
let $f_{*k}$ denote the induced map on $H_k(X;\Q)$. Note that $X_f$
fibers over
the torus $T^2$ with fiber $X$.  If  $\mbox{det }(I-f_{*1})=\pm 1$ then
there is a well-defined section class
$T$ (Corollary \ref{2.Twelldefined}).  Our first main result computes the
Gromov invariants of the multiples
of this section class.

\begin{theorem}\label{introThm1}
 If  $\mbox{det }(I-f_{*1})=\pm 1$,  the partial Gromov series of $X_f$ for
the section class $T$ is given by
the Lefschetz zeta function of $f$ in the variable $t=t_T$:
\bear
Gr^T(X_f)\ =\ \zeta_f(t)\ =\ \frac{\ma\prod_{k\  odd}\ \det(I-tf_{*k})}
{\ma\prod_{k\  even}\
\det(I-tf_{*k})}.
\label{zetafncthmIntro}
\eear
\end{theorem}

The proof occupies the middle part of this paper.  Section 3 reviews  the
zeta function, including  the
well-known second equality in (\ref{zetafncthmIntro}). We also  hint
at the first equality in
(\ref{zetafncthmIntro}) by showing  that the zeta function has a
factorization that parallels the factorization
of Taubes' Gromov invariant.  The difficult analysis is done in section
4, where we use a degeneration and gluing
argument to relate the Ruan-Tian invariants  to the zeta function.

\bigskip

 We can elaborate on this construction  by fiber summing $X_f$ with other
manifolds.  Suppose that $(Z,\w)$ is a
symplectic 4-manifold and  $F\subset Z$ is  a symplectically embedded
torus with trivial normal bundle.  Since $T\subset X_f$ is represented by a
symplectically embedded torus, we can identify
$T$ with $F\subset Z$ and form the fiber sum
\begin{equation}
Z(f)=Z\ma\#_{F=T} X_f.
\label{intro.fibersum}
\end{equation}
When this construction is done carefully, $Z(f)$ is a symplectic manifold
([G], [MW]).
\begin{theorem} The partial Gromov series of $Z(f)$ for the fiber class $F$ is
\bear
Gr^F(Z(f))\ =\ Gr^F(Z)\cdot \zeta_f(t_F)\cdot (1-t_F)^2.
\label{sumthmIntro}
\eear
\end{theorem}

\medskip

The righthand side of (\ref{sumthmIntro}) depends on $f$ only through its
induced  homology map.  It is
therefore easy to produce explicit examples.  In the last two sections we
give such examples in dimensions
four and six.

When  $X_f$ is a four-manifold,  a wealth of  examples arise from
knots.  Associated to each fibered knot $K$ in $S^3$ is a Riemann surface
$\Si$ and a monodromy diffeomorphism $f_K$ of $\Si$.  Taking $f=f_K$ gives
symplectic 4-manifolds $X_K$ of the homology type of  $S^2\ti T^2$ with
\bear\label{GrXK}
Gr(X_K)= \frac{A_K(t_T)}{(1-t_T)^2}
\eear
where $A_K(t)=\det(I-tf_{*1})$ is the Alexander polynomial of $K$ and $T$
is the section class.

More interesting examples are constructed by fiber summing with elliptic
surfaces.  Let $E(n)$ be the
simply-connected minimal elliptic surface with fiber $F$ and canonical
divisor $\kappa=(n-2)F$. Thus $E(1)$
is the rational elliptic surface and $K3=E(2)$; we will  also
write $E(0)=S^2\ti T^2$, even though this
is not simply-connected.   Forming  the fiber sum as in
(\ref{intro.fibersum}), we obtain a manifold
$$
E(n,K)\ =\ E(n)\#_{F=T} X_K.
$$
that  is symplectic and
homeomorphic to $E(n)$. In fact, for fibered knots $K,\; K'$ of the same
genus there is a homeomorphism between $E(n,K)$ and $E(n,K')$
preserving the periods of $\w$ and the canonical class $\kappa$.  These
manifolds are, however, distinguished by their Gromov invariants.
In this case we can compute the full series.

\begin{theorem}
\label{introThm2}
For $n\neq 1$, the Gromov invariant $E(n,K)$ is
\begin{equation}\label{introThm2.1}
Gr(E(n,K))\ =\ A_K(t_F)\ (1-t_F)^{n-2}
\end{equation}
and for $n>1$ this is also the Seiberg-Witten series of $E(n,K)$.
\end{theorem}
(See section 5 for comments on $E(0,K)$ and $E(1,K)$.) Thus fibered knots
with distinct  Alexander polynomials
give rise to symplectic manifolds $E(n,K)$ which are
homeomorphic but not diffeomorphic. In particular, there are infinitely
many distinct symplectic 4-manifolds
homeomorphic to $E(n)$.

Moving to dimension six, we can consider the symplectic manifolds
$E(n,K)\ti S^2$.  Using surgery
theory, one can show that for $n\geq 1$ fixed, the  $E(n,K)\ti S^2$
are diffeomorphic.  In fact, for knots of the same genus
these manifolds  have the same classical invariants  (Lemma
\ref{lemma7.1}).  On the other hand, in
section 6 we compute their Gromov
series and obtain the following.
\begin{theorem}
\label{introThm3}
For each $n\geq 1$, the smooth 6-manifolds $E(n)\ti S^2$ admit infinitely
many distinct symplectic structures,
all with the same classical invariants. In particular, this is true for
$K3\ti S^2$.
\end{theorem}
Here, `distinct' means that their symplectic forms are not deformation
equivalent. This simplifies and extends a result of
Ruan and Tian ([RT], Proposition 5.5). It is a striking example of the
``stabilization phenomenon'' suggested by Donaldson, namely that
homeomorphic but non-diffeomorphic symplectic 4-manifolds should give rise
to diffeomorphic but  deformation inequivalent symplectic 6-manifolds.

\bigskip

R. Fintushel and R. Stern have recently used knot theory and Seiberg-Witten
theory to prove a result
equivalent to Theorem \ref{introThm2}  above ([FS]).  This article arose
from our efforts to fit their results
into  a purely symplectic context and extend it to higher dimensions. We
thank both Rons for their help and
encouragement.
\bigskip

\bigskip

\setcounter{equation}{0}
\section{\bf Gromov invariants and Symplectic Sums}
\bigskip

The symplectic invariants we will use are  combinations of the invariants
defined by Ruan and Tian [RT].
To define them, we will briefly recall the construction of the Ruan-Tian
invariants,  specialize to
 ``stabilized degree 0'' RT invariants, and then explain how to
assemble these into the single symplectic invariant that we will call the
`Gromov series'.

\medskip

Let $(X,\w)$ be  a closed symplectic manifold of dimension $2n$. Choose a
compatible almost complex structure $J$;
this defines the canonical class $\kappa\in H^2(X)$ of $(X,\w)$. Given a
genus $g$ Riemann surface $\Si$ with
$k$ marked points $x_1,\dots,x_k$ and a homology class $A\in H_2(X,\Z)$,
consider pairs
$(j,f)$ where $j$ is a complex
structure on
$\Si$ and
$f$ is a map $\Si\to X$  satisfying
$$
\overline{\partial}_J f =\nu
$$
where $\overline{\partial}_J = df\circ j -J\circ df$ and where $\nu$ is a
fixed, appropriately defined  1-form.
The moduli space of such data $(j,x_1,\dots,x_k,f)$ has dimension $d$, where
\begin{equation}
\label{dimM}
d=2(n-3)(1-g)-2\kappa\cdot A +2k.
\end{equation}
We can reduce this dimension by requiring that $f$ take the marked points
into fixed constraint surfaces
$\al_i\subset X$. For appropriately chosen constraints, the moduli space
will be reduced to finitely many
points.  Counting with orientation, we get an invariant
$$
RT_{A,g,k}(\al_1,\dots,\al_k)
$$
that counts the total number of (perturbed) holomorphic
genus $g$ curves with homology class $A$ passing through the  constraints
$\al_1,\dots,\al_k$ (for generic almost
complex structure $J$ and  perturbation $\nu$).  This count depends only on the
deformation class of the symplectic structure ([RT] and [LT]).

\medskip

We will be interested only in the invariants for {\it unconstrained
curves}, i.e. for those pairs $(A,g)$  whose
$k=0$ moduli space has dimension 0.  Thus we define the ``degree zero''
invariants by
\begin{equation}
\label{1.defdegzero}
RT^0(A)\ =\ \sum_g RT_{A,g,0}
\end{equation}
where the sum is over all $g$ such that $\kappa\cdot A=(n-3)(1-g)$ (cf.
(\ref{dimM})).  Note that
there is at most one such $g$ when
$n\neq 3$, while for $n=3$ there is none unless $\kappa\cdot A=0$.

Definition (\ref{1.defdegzero}) must be clarified for the cases  $g=0$ and
$1$. In [RT], the invariants
$RT_{A,g,k}$ are defined  only for the ``stable range'' $2g+k\geq 3$.
However, the
definition can be extended to the unstable cases by choosing a class
$\beta\in H_{2n-2}(X)$ with $A\cdot \beta\ne 0$ and setting
\begin{equation}
RT_{A,1,0}\ =\ \frac{1}{A\cdot \beta}\, RT_{A,1,1}(\beta)
\qquad\mbox{and}\qquad
RT_{A,0,0}\ =\ \frac{1}{(A\cdot \beta)^3}\,
RT_{A,0,3}(\beta, \beta, \beta)
\label{unstable1}
\end{equation}
These invariants are independent of $\beta$ and still count perturbed
holomorphic curves (cf [IP1]).
\medskip

It is convenient to combine the  degree zero invariants $RT^0(A)$ into a
single quantity
associated with $X$.  To do
this, we introduce formal symbols $t_A$ for $A\in H_2(X;\Z)$ with relations
$t_{A+B}=t_At_B$ and construct a generating function. Our generating
function involves the M\"{o}bius $\mu$, which is defined by  $\mu(1)=1$
and for $m>1$
\[ \mu(m)=\l\{\begin{array}{cl} (-1)^k & \mbox{  if $m$ is a product
of $k$ distinct primes }\\
0&\mbox{ otherwise }
\end{array}\r.\]

\begin{defn} The {\em degree zero Gromov series } is
\bear\label{defgr}
Gr^0(X)= \mbox{\em exp}\left[\sum_{A}  RT^0(A) \ \phi(t_A)\right]
\eear
where $\phi(t)=\sum_{k=1}^\infty \mu(k)t^k$ is the power series with
coefficients given by  the M\"{o}bius function.  For $A\in H_2(X)$, the
 Gromov invariant $Gr^0_X(A)$  is defined by the expansion $Gr^0(X)=\sum
Gr^0_X(A)\,t_A$.
\end{defn}
 Of course, different weighting functions $\phi$ give different symplectic
invariants.  For the simple
choice $\phi(t)=t$, $ Gr^0_X(A)$ counts the
number of ways we can represent the homology class $A$
as the image of a perturbed degree zero holomorphic map whose domain is a
disjoint union of (stable) Riemann surfaces.  While it initially seems
awkward to introduce the M\"{o}bius
function  into (\ref{defgr}), this choice eliminates some ``overcounting'',
giving a simpler generating function which contains the same information.

\medskip

\rem\ \ One can arrive at (\ref{defgr}) as follows.  Fix  a  function $F$,
and assign a factor
$F(t_A)$ to each curve that contributes $+1$ to the count $RT^0(A)$, and a
factor $1/F(t_A)$ to each curve that
contributes $-1$, obtaining
\bear\label{predefgr}
Gr^0(X)=\ma\prod_{A\in H_2(X)}F(t_A)^{RT^0(A)},
\eear
Again, the natural choice  $F(t)=e^t$ overcounts; it is more efficient to
choose $F$ so that
\bear\label{prodF=t}
\prod_{k=1}^\infty F(t^k)=e^t.
\eear
This holds for
\bear\label{MoebiusF}
F(t)=\exp\l(\ma\sum_{m\ge 1} \mu(m)t^m \r),
\eear
as can be verified   by writing $\ell=mk$ and using the
basic fact that  $\ma\sum_{m|\ell}\mu(m)=0$ unless $\ell=1$, in which case
the sum equals to 1.

\bigskip

We will usually omit the superscript from $RT^0(A)$, $Gr^0(X)$, and
$Gr_X^0(A)$, it being understood
that we are always computing degree 0 invariants.  We will often be
interested in the invariants
$Gr_X^0(mA)$ for multiples of a given class $A$; it is then convenient to
work with the {\it partial Gromov series}
\bear
\label{defpartialgr}
Gr^A(X)= \ma\prod_{m\ge 1}F(t_A^m)^{RT^0(mA)}
\eear
and then
\best
Gr(X)=\ \ma\prod_{A\; primitive} Gr^A(X).
\eest

\bigskip

The Gromov invariant (\ref{defgr}) has several nice properties.  In
particular, it  behaves well under the
`symplectic normal  sum' construction described by Gompf [G] and
McCarthy-Wolfson [MW].  A general formula for the
invariants of the sum is given in [IP2].  Here we will state and use the
formula only in the special case of an
ordinary fiber sum.

Fix a  $2n$-dimensional symplectic manifold $X$, not necessarily connected.
By a {\it symplectic fiber} in
$X$ we mean  a codimension 2 symplectic submanifold with trivial normal
bundle.  Suppose that $F_1$ and $F_2$ are
 disjoint symplectic fibers in $X$ which are symplectomorphic. Following
[G] and [MW], we can  construct an
orientation reversing symplectomorphism
\bear
\phi:N' F_1\ra N'F_2
\eear
where $N'F=N_{\ep}F\setminus F$ denotes a $2\ep$-tubular neighborhood (for
some metric) with its core removed.
Removing $\ep$-neighborhoods of $F_1$ and $F_2$ and identifying the
boundaries via $\phi$, yields the
 symplectic fiber   sum  $\#_\phi X$.  This is a symplectic manifold, and the
  deformation type of its symplectic structure depends on the deformation
class of $\phi$. We will usually fix an deformation class of $\phi$ and
denote $\#_\phi X$ by ${\#}_{F_1=F_2} X$, or even $\#_{F} X$. When $X$ is
the disjoint union of manifolds
$X_1$ and $X_2$ and $F_i\subset X_i$ as above, we denote the corresponding
fiber sum by
$$
X_1 \#_{F_1=F_2}  X_2 \qquad\mbox{or}\qquad  X_1 \#_{F}  X_2.
$$
Note that in the sum, the homology classes $[F_1]$ and $[F_2]$ become a
single `fiber class' $[F]$.

\bigskip

\begin{theorem}(cf. [IP2])  Let $X_1$ and $X_2$ be closed symplectic
4-manifolds containing symplectic fibers
$F_i\subset X_i$ as above.  Then the partial Gromov invariant for the fiber
class $[F]$ is
\bear
Gr^{[F]}_{X_1\#_F X_2}=(1-t_{F})^2 \cdot Gr^{[F_1]}_{X_1}\cdot
Gr^{[F_2]}_{X_2}.
\label{1.gluinglemma}
\eear
\end{theorem}

\medskip\bigskip

As an application, let us compute the Gromov series of $E(n)$. By
an observation of Donaldson [D], a generic K\"{a}hler structure $J$ on
$E(2)=K3$ admits no holomorphic curves whatsoever. By a limiting argument,
this implies that $RT(A)= 0$
for each non-trivial class $A\in H_2(K3)$. Specifically, if $RT(A)\neq
0$,  we could choose a   sequence of
generic $(J,\nu)$ converging to some $(J_0,0)$ with $J_0$ generic and
K\"{a}hler and a sequence of
$(J,\nu)$-holomorphic A-curves.  These curves would limit to  a bubble tree
of $J_0$-holomorphic curves (cf.
[PW], [P], and [RT]), contradicting Donaldson's observation. Thus
\begin{equation}\label{1.K3}
Gr(K3)=1.
\end{equation}
Now, as symplectic manifolds
\begin{equation}
E(n)\#_F E(1)=E(n+1)
\label{1.EmEnEmn}
\end{equation}
where this fiber sum glues a fiber of $E(n)$ to a fiber of $E(1)$. Taking
$n=1$ and applying (\ref{1.gluinglemma}) yields $Gr^F(E(1))=1/(1-t)$;
putting this into
(\ref{1.EmEnEmn}) and (\ref{1.gluinglemma}) inductively gives
\begin{equation}
Gr^F(E(n))\ =\ (1-t)^{n-2}
\label{1.GTEn}
\end{equation}
where the variable $t$ corresponds to the fiber class.  In fact, we will
show in Lemma \ref{fullGr=Grt} that
for $n\geq 2$, this formula gives the {\it full} Gromov series of $E(n)$.

\bigskip

\bigskip

\setcounter{equation}{0}
\section{\bf The Symplectic Mapping Torus}
\bigskip

Fix a closed symplectic manifold $(X,\w)$.  Each symplectic diffeomorphism
$f:X\to X$ induces a
symplectomorphism $F$ of $X\times \R\times S^1$ (with the product
symplectic structure) by
$F(x,s,\theta)=(f(x),s+1,\theta)$.  The quotient
\begin{equation}
X_f = {X\times \R\times S^1}/\Z
\label{3.defXf}
\end{equation}
is a closed symplectic manifold --- the {\it symplectic mapping torus} of
$f$.   As a differential manifold,
$X_f$ is the product $M_f\times S^1$, where $M_f$ is the usual mapping
torus. Note that $X_f$ fibers over
the torus $T^2$ with fiber $X$, and that  the $S^1$ factor gives a
 symplectic circle action on $X_f$.

We can also construct $X_f$ as a symplectic fiber sum. Start with
$S^2\ti X$ and fix two points $a,\; b\in S^2$. Let $X_a=\{a\}\ti X$ and
$X_b=\{ b\}\ti X$. Extend  $f:X_a\ra X_b$ to a symplectic map between the
$\ep$-tubular neighborhoods of $X_a$ and $X_b$, and form the symplectic
fiber sum.  The result is $X_f$:
\bear
\underset{f(X_a)=X_b}{\Huge{\#}} S^2\ti X\ =\ X_f.
\label{3.2}
\eear
To see this equivalence, put a linear circle action on $S^2\subset\R^3$ and
take $X_a$ and $X_b$ to be the
fibers at the poles, and let $S^2_{\ep}$ be the sphere with $\ep$-disks
around the poles removed.  While
the  usual symplectic structures on
$S^2_{\ep}$ and $[0,1]\ti S^1$ are not symplectomorphic, they are isotopic,
and hence (\ref{3.defXf}) and
(\ref{3.2}) are isotopic symplectic manifolds.

\bigskip

\begin{ex}{\rm Let $f$ be the symplectomorphism of the standard
2-torus induced by the matrix
$$
A=\left(\begin{array}{cc}1&1\\0&1\end{array}\right)
$$
acting on the torus $\R^2/\Z^2$.  Then $\mbox{cok}
(I-f_{*1})=\mbox{cok} (I-A)=\Z$, so by Lemma \ref{piX_f} below, $X_f$ has
first betti number $b_1=3$.
This is Thurston's famous example [Th] of a symplectic manifold with no
K\"{a}hler structure (compact
K\"{a}hler manifolds have
$b_1$ even).}
\end{ex}

\bigskip

We can always assume that $f$ has a fixed point $p$ (if not, replace $f$ by
 its composition with a
Hamiltonian flow taking $f(p)$ back to
$p$; this changes $X_f$ to an isotopic
symplectic manifold).  The differential $df_p$ at the fixed point lies in
the connected group
$\mbox{Sp}(T_pX)$.  Consequently, $\{p\}\ti \R\ti S^1/\Z$ defines a section
$T_p$
\begin{equation}
\begin{array}{cc}
X_f  \\
\downarrow \uparrow\vcenter{\rlap{$T_p$}}\\
T^2 \\
\end{array}
\label{3.sectionmap}
\end{equation}
whose image is a torus with trivial normal bundle.  Let $T=[T_p]$
denote this `section
class' in $H_2(X_f;\Z)$ (in general, $t$ depends on the choice of the fixed
point $p$).

\bigskip

It is straightforward to calculate the fundamental group of $X_f=M_f\ti S^1$.

\begin{lemma}
 Let $f_{\#}$ and $f_*$ be, respectively, the induced maps on
$\pi_1(X)$ and $H_1(X,\Z)$.
Then
\begin{enumerate}
\item[(a)] $\pi_1X_f = G\times \Z$ where $G=\{\langle \pi_1X,\tau\rangle\ |\
\tau^{-1}x \tau=f_{\#}x\quad \ \forall
x\in \pi_1X \}$, and\
\item[(b)] $H_1(X_f;\Z)= \mbox{cok} (I-f_{*1})\oplus\Z\oplus\Z$.
\end{enumerate}
\label{piX_f}
\end{lemma}
\pf A fixed point $p$ in $X$ determines  a path $\tilde{\tau}=\{p\}\ti
[0,1]$ in $X\ti [0,1]$ that projects to
a loop $\tau$ in $M_f$. Then $\pi_1M_f$ is generated by $\pi_1 X$ and
$\tau$ and the only relations are those
obtained by conjugating by $\tau$; namely,
for each $x\in\pi_1(X;p)$ the loops $\tilde{\tau}(f_{\#}\times
\{1\})\tilde{\tau}^{-1}$ and $x\times \{0\}$ are
homotopic rel basepoint in $X\times \R$, so $\tau (f_{\#}x) \tau^{-1}=x$ in
$\pi_1M_f$.  This gives (a).

Abelianizing, we have $H_1(X_f)= K\oplus\Z\oplus\Z$ where $K$ is the quotient
of $H_1(X;\Z)$ by the relation $x=f_{*1}x$. This gives (b). \qed

\medskip

The proof of Lemma \ref{piX_f} also implies the following simplification.
\begin{cor}
If $\mbox{det }(I-f_{*1})=\pm 1$ then the section class $T=[T_p]\in
H_2(X_f;\Z)$ is well-defined, independent of the fixed point $p$.
\label{2.Twelldefined}
\end{cor}

\bigskip

We now focus on the case where $X$ is a Riemann surface of genus $g$. Each
orientation-preserving diffeomorphism $f$ of $X$ preserves the cohomology
class of the symplectic form.  By
Moser's Theorem [M] we can assume (after an isotopy of $f$) that $f$ is a
{\it symplectic} diffeomorphism,
and hence determines  a symplectic mapping torus $X_{f}$. Write $f_{*1}$ for
the induced map on $H_1(X,\Z)$ and set
\bear\label{A=det(I-tf)}
A(t)=\mbox{det}\; (I-tf_{*1})
\eear
\begin{lemma}
\label{l2.4}
(a) $A(t)$ is a monic polynomial with integer coefficients with
$A(t^{-1})=t^{-2g}A(t)$.
\indent\hskip.68in  (b) $X_f$ is a homology $S^2\ti T^2$ if and only if
 $A(1)=\pm 1$.
\end{lemma}
\pf The induced map $f_{*1}$ is symplectic by Poincar\'{e} duality. After
fixing a basis of $H_1(X,\Z)$ the matrix $A=f_{*1}$ lies in
$Sp(2g,\Z)$, so $\det A=1$.  Then the first two properties listed in (a) are
obvious. For the third, note that $A(t)$ is invariant under
conjugation so it suffices to verify $A(t^{-1})=t^{-2g}A(t)$ for matrices in
the maximal torus of $Sp(2g,\cx)$.  But that is clear since such matrices
 have the form
$\mbox{diag}(\lambda_1,\lambda_1^{-1},\dots,\lambda_g,\lambda_g^{-1})$.

Finally, (b) holds because $X_f=M_f\ti S^1$ and, by Lemma \ref{piX_f}(b) and
Poincar\'{e} duality,  $M_f$ is a homology $S^2\ti S^1$ if and only if
$I-f_{*1}$
is invertible over $\Z$. \qed

\bigskip

Fibered knots provide one source of diffeomorphisms of Riemann surfaces.
Recall that a knot $K\subset S^3$ is
{\it genus $g$ fibered}   if there is an oriented fibration
$\pi_K:S^3\setminus K\to S^1$ whose fiber $X_0$
is the 2-manifold of genus $g$ with one point removed. After a 0-surgery, this
becomes an $X$-fibration over $S^1$. This new fibration is the mapping torus
$M_{f_K}$ of the monodromy diffeomorphism $f_K:X\ra X$ along the knot,
$M_{f_K}$ is a homology $S^2\ti S^1$ and
\bear\label{Mf-tau}
S^3\setminus K=M_{f_K} \setminus \{\tau\}
\eear
Moreover, the Alexander polynomial of $K$ is the characteristic polynomial
of the map induced by $f_K$ on $H_1(X)$ exactly as in (\ref{A=det(I-tf)}).
\bigskip

\begin{ex}{\rm The monodromy matrices $f_*$ of the trefoil and figure 8
knots are, respectively,
$$
\left(\begin{array}{cc}1&1\\-1&0\end{array}\right)
\hskip1in
\left(\begin{array}{cc}2&1\\1&1\end{array}\right).
$$
and the Alexander polynomials are $t^2-t-1$ and $t^2-3t+1$.
}\end{ex}

\bigskip

Looking the other way, Burde [B] proved that any polynomial $A$ as in
 Lemma \ref{l2.4}(a) with $A(1)=\pm 1$ arises as the Alexander polynomial
 of a knot in $S^3$. Summarizing,
\begin{lemma}\label{l2.5} For each genus $g$ fibered knot in a homology
3-sphere, there is a symplectic 4-manifold
\best
X_K=X_{f_K}
\eest
This is a homology $S^2\ti T^2$ and $\mbox{det}(I-t(f_K)_{*1})=A_K(t)$ is a
monic
polynomial with integer coefficients satisfying $A_K(1)=\pm 1$ and
$A_K(t^{-1})=t^{-2g}A_K(t)$. Conversely, each such $A_K(t)$ arises as
$\mbox{det} (I-t(f_K)_{*1})$ for some such $K$.
\end{lemma}


\bigskip
\setcounter{equation}{0}
\section{\bf Factorizing the Lefschetz Zeta Function}
\bigskip

The Lefschetz number $L(f)$ of $f:X\to X$ is the homological intersection
of the graph of $f$ with the diagonal in $X\times X$. When this
intersection is transversal
\bear\label{Lffixedpt}
L(f)=\sum_{f(x)=x}L(f,x) \quad \mbox{ where } \quad
L(f,x)=\mbox{ sgn det}(\mbox{df}_x-I)
\eear
The Lefschetz trace formula asserts that
\bear
L(f)=\sum_{k=0}^{dim\ X} (-1)^k\ \tr(f_{*k})
\label{Lftraceform}
\eear
where $f_{*k}$ is the induced endomorphism of $H_k(X;\Q)$.  Often, as in
dynamics, one considers the
Lefschetz numbers $L(f^n)$ of the iterates of $f$.  These are neatly
encoded in the
{\it Lefschetz zeta function}
\bear
\zeta_f(t) \ =\ \exp\left({\sum_{n=1}^{\infty}\frac{t^n}{n} L(f^n)}\right).
\label{zetafnc}
\eear
Substituting in (\ref{Lftraceform}), one sees that this sum converges to a
continuous function for small $t$.

Starting form formulas (\ref{Lffixedpt})-- (\ref{zetafnc}), we will obtain
two formulas for the zeta function.  The first --- the homological
expression (\ref{zeta=dets}) --- is well-known.  The second expresses
$\zeta_f$ as a product of generating functions.  The arguments required for
this are completely elementary, but the resulting ``factorization formula''
does not seem to be in the
literature.

\begin{lemma}
$\zeta_f$ extends to a rational function of $t\in\cx$. In fact,
\bear
\zeta_f(t)\ =\ \frac{\ma \prod_{k\  odd}\ \det(I-tf_{*k})}
{\ma \prod_{k\  even}\
\det(I-tf_{*k})}.
\label{zeta=dets}
\eear
\label{zetaformulalemma}
\end{lemma}
\pf For any $n\times n$ complex matrix $A$ and  $t<<|A|$ we have the identity
\bear
\exp\left({-\sum_{n=1}^{\infty} \frac{t^n}{n} \tr A^n}\right)\ =\ \det(I-tA).
\label{zetafnc2}
\eear
In the simplest case, when $A$ is a $1\times 1$ matrix,  (\ref{zetafnc2})
is the Taylor series of
$\ln(1-tA)$.  The general case follows from three observations: (i)  the
functions on both sides of
(\ref{zetafnc2}) satisfy
$f(A\oplus B)=f(A)f(B)$, so (\ref{zetafnc2}) holds for diagonal matrices,
(ii) the functions on both sides
are invariant under conjugation, so (\ref{zetafnc2}) holds for
diagonalizable matrices, and (iii) by Jordon
normal form any matrix is a limit of diagonalizable matrices, and both
sides of (\ref{zetafnc2}) are
continuous.

Equation (\ref{zeta=dets}) is easily obtained from
(\ref{Lftraceform}), (\ref{zetafnc}) and  (\ref{zetafnc2}). \qed
\medskip\bigskip

A second expression for $\zeta_f$ arises by substituting (\ref{Lffixedpt})
into (\ref{zetafnc}).  At a fixed point $p$ of a symplectic diffeomorphism
$f$, the differential
$df_p:T_pX\to T_pX$ is real and
symplectic, and contributes $L(f,x)=\mbox{sgn det }(df_p-I)=\pm 1$ to the
Lefschetz number.
After examining the structure of the symplectic group, we will be able to
compute these signs for the powers of $f$.
This will yield an expression for the Lefschetz zeta function as a
product of generating functions
associated to the finite orbits of $f$.

Let $A\in Sp(2n,\R)$, i.e. $A^tJA=J$ where $J$ is the $2n\ti 2n$ matrix
$\left(\begin{smallmatrix}0&-I\\ I&0\end{smallmatrix}\right)$. When $n=1$
such $A$ are either
elliptic or hyperbolic type --- in some basis they have one of the  forms
$$
E=\begin{pmatrix}\cos\theta & \sin\theta\\ -\sin\theta &\cos\theta\end{pmatrix}
\hskip1in
H=\begin{pmatrix}\lambda & 0 \\ 0 & \lambda^{-1}\end{pmatrix}\qquad\lambda
\in\R.
$$
In general,  $A\in Sp(2n,\R)$ can be put into one of three block forms (cf.
[R2]):
\bear
\begin{pmatrix}E&0\\ 0&S\end{pmatrix}, \qquad
\begin{pmatrix}H&0\\ 0&S\end{pmatrix}, \qquad
\begin{pmatrix}E&* &*\\ 0&E'&0\\ 0& *&S\end{pmatrix}.
\label{3.blocks}
\eear
In the first two cases $E$ and $H$ are as above and $S\in Sp(2n-2,\R)$,
while  in the last case $E,E'$ are
elliptic with eigenvalues
$\lambda,\overline{\lambda},\lambda^{-1},\overline{\lambda}^{-1}$ and $S\in
Sp(2n-4,\R)$.  Thus (\ref{3.blocks}) inductively defines a `normal form'
for symplectic matrices.

\medskip

The matrices that have 1 as an eigenvalue define a codimension-one
algebraic subset $W_1$ of $Sp(2n,\R)$.
Each $A\notin W_1$ has a certain number $R(A)$  of real eigenvalues,
$R^+(A)$ of which are positive.  It is
easy to see that $Sp(2n,\R)\setminus W_1$  has three open components:
\begin{equation}
\left\{
\begin{array}{ll}
\mbox{type } E:   &  \mbox{those with $R(A)$ even,} \\
\mbox{type } H:   &  \mbox{those with $R(A)$ odd but $R^+(A)$ even,}\\
\mbox{type } H':  &   \mbox{those with $R(A)$  and  $R^+(A)$ both odd.}
\end{array}
\right.
\label{3types}
\end{equation}
In the simplest case, a $2\ti 2$ matrix $A\in Sp(2,\R)$ has type $E$, $H$,
or $H'$ when it is, respectively,
elliptic, hyperbolic with positive eigenvalues, or hyperbolic with negative
eigenvalues.

 We can also consider the  walls $W_m=\{A\ |\ \mbox{det
}(A^m-I)=0\}$ formed by those $A$ with some
eigenvalue equal to a $m^{\mbox{th}}$ root of unity.  Then the complement
of $W=\cup_{m\neq 0} W_m$ is a Baire
subset of $Sp(2n,\R)$.

Define the sign of $A\in Sp(2n,\R)\setminus W_1$ by
$$
\mbox{sgn }(A) = \mbox{sgn det }(A-I).
$$

\begin{lemma}
If $A\in Sp(2n,\R)\setminus W$, then for each non-zero $m\in \Z$,
$$
{\mbox{ sgn }} (A^m)\ =\ \left\{
\begin{array}{cl}
1  & \mbox{ if $A$ has type $E$,}\\
-1  & \mbox{ if $A$  has type $H$,}\\
-(-1)^m  & \mbox{ if $A$  has type $H'$.}\\
\end{array}
\right.
$$
\label{spLemmma1}
\end{lemma}
\pf  From the normal form (\ref{3.blocks}) we have $\mbox{sgn } A^m=\prod
\mbox{sgn det }
(B^m_i-I)$ where $B_i$ are  $2\ti
2$ blocks along the diagonal.  Each elliptic block has eigenvalues
$\lambda,\overline{\lambda}$, so has
$\mbox{det } (B^m-I)=|\lambda^m-1|^2>0$ because $\lambda$ is not a root of
unity.  Similarly, for each
hyperbolic block  the sign of $\mbox{det }(B^m-I) =
(\lambda^m-1)(\lambda^{-m}-1)= -|\lambda^m-1|^2/\lambda^m$
is $-(\mbox{sgn }\lambda)^m.$ The Lemma follows using the definitions
(\ref{3types}).
\qed

\medskip

To apply these sign counts, we use the following transversality result,
which follows from a result of R.C.
Robinson ([R1] Theorem 1Bi) and the fact that the group of symplectic
diffeomorphisms is locally connected.

\begin{lemma}
Each symplectic diffeomorphism $f_0:X\to X$ is isotopic through symplectic
diffeomorphisms to a symplectic
diffeomorphism $f$ for which $\mbox{graph}(f^n)$ is transverse to the
diagonal for all $n\geq 1$.
\label{transversality}
\end{lemma}

Since a symplectic isotopy changes neither the Lefschetz numbers nor the
isotopy class of the symplectic
manifold $X_f$, we may henceforth assume that the powers of $f$ are
transverse as in Lemma
 \ref{transversality}.   At a fixed point $p$ of $f$, this transversality
condition says that
$\mbox{det }(df^n_p-I)\neq 0$.  Thus the signs of $df_p$ and $df^n_p$ are
related
as in Lemma \ref{spLemmma1}.

\medskip

Now each orbit of $f$ of minimal period $k$ has exactly $k$ points, all of
the same type (\ref{3types}).  Set
$$
\left\{
\begin{array}{ccl}
e_k & = &  \mbox{the number of type $E$ orbits of $f$ of minimal period $k$,} \\
h_k & = &  \mbox{the number of type $H$ orbits of $f$ of minimal period $k$,} \\
h'_k & = &  \mbox{the number of type $H'$ orbits of $f$ of minimal period $k$.}
\end{array}
\right.
$$
Then, because each fixed point of $f^n$  has a minimal period $k|n$, we have
$$
L(f^n)\ =\ \sum_{k|n} k\left[ e_k-h_k-(-1)^{n/k} h'_k\right].
$$
Substituting this into the Lefschetz zeta function (\ref{zetafnc}), writing
$n=kl$, and rearranging the sum
gives
$$
\zeta_f(t)\ =\ \mbox{exp}\l[ \sum_{k=0}^{\infty} \sum_{l=0}^{\infty}\left[
e_k-h_k-(-1)^{n/k}
h'_k\right] \frac{t^{kl}}{l}\r].
$$
Using the Taylor series  $\log(1-t^k)=-\sum t^{kl}/l$, this reduces to
\begin{equation}
\zeta_f(t)\ =\ \prod_k
\left(\frac{1}{1-t^k}\right)^{e_k-h_k}\left(1+t^k\right)^{h'_k}
\label{taubesprod1}
\end{equation}
In this we see the three ``generating functions''
\begin{equation}
f_E(t)=\frac{1}{1-t} \qquad f_H(t)=1-t, \qquad f_{H'}(t)=1+t.
\label{lefgeneratingfncs}
\end{equation}
Thus if we let ${\cal O}$ be the set of finite-period orbits of $f$, and
assign to each $\tau\in{\cal O}$ of
order $k=k(\tau)$ one of the generating function
$f_{\tau}(t)=f_E(t^k),f_H(t^k)$ or $f_{H'}(t^k)$ according to
the type of $\tau$, we obtain
\begin{equation}
\zeta_f(t)\ =\ \prod_{{\tau}\in{\cal O}} f_\tau(t^{k({\tau})}).
\label{taubesprod2}
\end{equation}

\medskip

\begin{rem}  Bifurcations in the dynamics of $f$ give relations amongst the
generating functions (\ref{lefgeneratingfncs}). In fact, for a
one-parameter family $\{f_s\}$ of
symplectic diffeomorphisms, two
types of stable bifurcations can occur (cf. [AM] \S 8.6--8.7).

\begin{enumerate}
\item[(i)]  At some $s_0$ a new pair of finite-order points is created, one
of type $E$ the other of type $H$,
both  with the same minimal  period $k$.  The invariance of the zeta
function (\ref{taubesprod2}) then corresponds to the relation
\begin{equation}
f_E(t^k)f_H(t^k)=1.
\label{4.5}
\end{equation}

\item[(ii)]  A more interesting bifurcation is called ``subtle division'':
at some $s_0$ a type $H'$ orbit of
minimal period $k$  disappears, and is replaced by a type $E$ orbit of
minimal period $k$ and  a type $H$
orbit of minimal period $2k$.  The corresponding  relation is
\bear
f_{H'}(t^k)=f_E(t^k)f_H(t^{2k}).
\label{4.6}
\eear
\end{enumerate}
\end{rem}

\medskip

The generating functions (\ref{lefgeneratingfncs}) are identical to three
of the generating functions Taubes
uses to define his Gromov invariants [T3]. In fact, Taubes derives his
generating functions using  exactly the relations (\ref{4.5}) and (\ref{4.6}).
  Apparently  there is a precise correspondence between the periodic orbits
of $f$ and the $J$-holomorphic tori in $X_f$, with the orbits of types $E$,
 $H$, and $H'$ corresponding to tori of Taubes type $(+,0)$, $(-,0)$, and
 $(+,1)$ respectively (cf. [T3]).  That is the content of our main result
Theorem \ref{introThm1}.  For the proof, however, we must look away from
the dynamics of $f$ and instead study
the holomorphic curves in $X_f$.

\bigskip


\bigskip
\setcounter{equation}{0}
\section{\bf The Gromov Series of $X_f$}
\bigskip

 In this section we will use
the Gluing Theorem of [IP2] to establish the following remarkable formula
linking the Gromov series of $X_f$ to the dynamics of the diffeomorphism $f$.

\begin{theorem} For dim $X_f=2n\neq 6$, the  partial Gromov Series of
$X_f$ with respect with $T$ is equal to the zeta function of $f$:
\bear
Gr^T_{X_f}\ =\ \zeta_f
\label{5.zetafncthm}
\eear
When $n=3$ this formula holds when the righthand side is the ``$g=1$
partial Gromov series'' defined by
(\ref{defpartialgr}) with $RT(kA)$ replaced by $RT_{kA,1,0}$, i.e. by the
genus one term in
(\ref{1.defdegzero}).
\label{theorem5.1}
\end{theorem}
This, together with the zeta function  formula (\ref{zeta=dets}), yields
Theorem \ref{zetafncthmIntro}  of the introduction.

\medskip

To start, note that the canonical bundle of $X_f$ is trivial along the
section defined in equation (\ref{3.sectionmap}), so $\kappa\cdot T=0$.  The
dimension formula (\ref{dimM}) then shows that,  when $n\neq 3$, the
only curves in $X_f$ that
contribute to the series (\ref{5.zetafncthm})  are tori; when $n=3$ we will
simply restrict attention to
tori.  To count these, we must stablize as in
(\ref{unstable1}). For this, we fix a  point $c\in T^2$,  let $X_c$
be the fiber of $X_f\to T^2$ over $c$, and consider the moduli space
\bear
\M_{mT}(T^2,X_f)
\label{5.3}
\eear
of perturbed holomorphic maps $\phi$ from the generic torus with a marked
point $x_0$ to $X_f$ which
represent $mT$ in homology and satisfy $\phi(x_0)\in X_c$.  For generic
$(J,\nu)$, this moduli space has
$$
RT_m = RT_{mT,1,1}([X_c])
$$
oriented points.  Referring to (\ref{defpartialgr}) and (\ref{unstable1}),
and noting that $X_c\cdot mT=m$,
\bear
Gr^T&=&\ma\prod_{m=1}^\infty \l[  F(t^m)\r]^{ RT_m/ m}
\label{5.2}
\eear
Thus we must calculate the invariant $RT_m$.

\medskip

The key ingredient is  a gluing theorem, which relates the moduli space
(\ref{5.3}) of curves in  the glued
manifold $X_f$ to a moduli space of curves in  the ``unglued'' manifold
$S^2\ti X$.  This result
(Proposition \ref{RT convergence}) is largely subsumed by the general
gluing theorem given in [IP2].  With that
in mind, we will outline the part of the argument that overlaps with [IP2],
and give complete details on
those points particular to the present case.

\medskip

 Recall that $X_f$, more
properly denoted $X_f^{\ep}$, is
constructed from $S^2\times X$ by fixing points $a,b\in S^2$, identifying
the fibers $X_a, X_b$ via $f$,   and
symplectically gluing the  $\ep$-tubular neighborhoods
 $N_a(\ep)$ and $N_b(\ep)$.  The degeneration formula describes what happens
  to holomorphic curves as  $\ep\to 0$ and explains how to compute the
  Gromov invariants from the limiting curves.

As $\ep\to 0$, the spaces $X_f^{\ep}$ converge to
$X_f^{0}=S^2\times X/\sim$,
 where $\sim$ is the identification of $X_a$ with $X_b$ via $f$.
This limiting space is singular, but is a manifold away from $X_a=X_b$.
 Similarly, as $\ep\to 0$, the curves in
$X_f^{\ep}$ representing a multiple  $mT$ of the section class
$T\in H_2(X_f^{\ep})$ converge to singular curves in  $X^0_f$ representing
$mT_0\in H_2(X_f^0)$.

 More precisely, let $f_\ep:(T^2,j_\ep)\to X$  be a sequence of
$(J,\nu)$-holomorphic maps
in $X_f^{\ep}$ that take a marked point $x_0$ into a fixed fiber $X_c$
(as above, the domains must be tori). As
$\ep\ra 0$, the sequence of complex structures
$j_\ep$ has a  subsequence converging to some $j_0$ in the Deligne-Mumford
compactification of the moduli space of complex structures on
$(T^2,x_0)$.  The Bubble Tree Convergence Theorem ([PW], [RT])
implies that $(f_\ep,j_\ep,x_0)$ has a subsequence that converges modulo
diffeomorphisms to a map
$f:B\ra X^0_f$ where $B$ is a bubble tree.

The domain of a bubble tree $B$ is a connected union of three
types of components (cf. Figure 4.1).  The original torus limits to a
``principal component''
$(T^2,j_{\infty}, x_0)$, which is either a  torus or a sphere with two
special points besides $x_0$.  There
can be finitely many trees of ``ordinary bubbles'' attached to the
principal component.  Finally, if the
principal component is singular, there can be  ``chains'' of bubble
spheres, each
with exactly two special points, inserted between the two special points of
the principal domain. The perturbation $\nu$ extends as 0
on all the bubbles.

\vskip.2in
\centerline{\psfig{figure=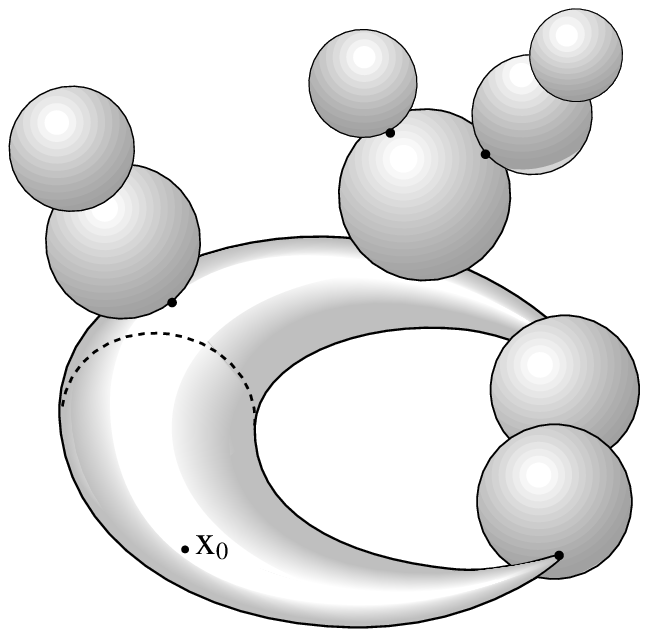,width=2.0in}}
\begin{center}{Figure  4.1 --- A general bubble tree domain}\end{center}
\vskip.2in

To simplify the exposition, we first assume that the limit has no
components entirely contained
in the singular set $X_a=X_b$ of $X_f^0$. Then the limiting map has a
special property: the double points of
the domain are exactly those points in the inverse image
of the singular set [IP2]. Remove all these double points and lift each
component $B_i$ of $B\setminus D$ to a map $B_i\ra  S^2 \ti X$. By the
Removable Singularity Theorem for
$J$-holomorphic curves [PW], these extend to $J$-holomorphic maps $\ov
B_i\ra S^2\ti X$. Thus
we may consider the limit as a map
\bear
\label{5.liftedBT}
\bigsqcup \ov B_i\to S^2\ti X
\eear
together with the constraint that the image intersects  $X_a$ and $X_b$ in
the same set
of points (counted with multiplicities).

\medskip

Since no homology class in $H_2(S^2\ti X)$ can intersect $X_a\cup X_b$
in exactly one point,  there are no
ordinary bubbles (each tree of ordinary bubbles has a top bubble
with only one double point).  It also means that the principal component is
 a sphere ($j=\infty$), since otherwise it is a torus with no double
points.   Thus the domain of the limit
map is a  ``necklace'' of spheres, one of them marked, and the others  having
exactly two double points each, one mapped  into $X_a$ and the other mapped
into $X_b$ (Figure 4.2).

\vskip.3in
\centerline{\psfig{figure=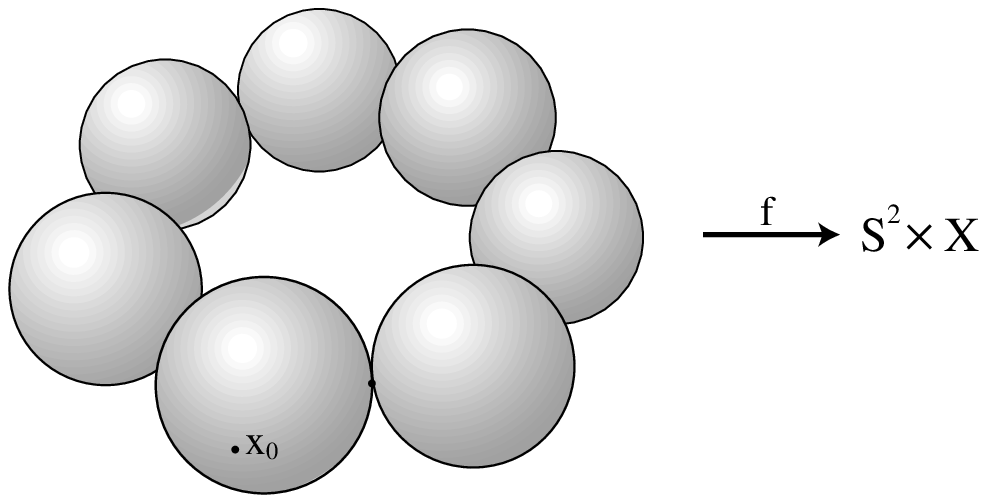,width=3.0in}}
\begin{center}{Figure  4.2 ---  A necklace map}\end{center}
\vskip.3in

Assume the the necklace has $k$ beads, the $i$-th bead representing  $A_i$
in homology.
The constraint that the image intersects  $X_a$ and $X_b$ in the same set
of points (counted with multiplicities) implies that each contact point of
the necklace with  $X_a$ and $X_b$ has multiplicity
$d=A_i\cdot X_a=A_i\cdot X_b$, {\em for the same $d$}.   Thus we can write
$A_i=ds+\alpha_i\in H_2(S^2\ti X)$
where $s=[S^2\ti \mbox{pt}]$ and $\alpha_i=[
\mbox{pt}\ti \alpha_i']$ for some $\alpha'_i\in H_2(X)$.  Furthermore, the total
lifted bubble tree represents an element of
$\pi_*^{-1}(mT_0)=ms+\mbox{ker}\,\pi_*$
in $H_2(S^2\ti X)$, where $\pi$ is the projection $S^2\ti X\ra X_f^0$.
Therefore we have
\bear\label{5.sumalpha=0}
m=kd \qquad \mbox{and}\qquad   \sum\alpha'_i\in \mbox{ker}\, \pi_*.
\eear

\medskip

Note that only one bead of the necklace is a stable curve. To make  the
gluing construction well-defined, we stabilize the other
$k-1$ beads as in (\ref{unstable1}), using $\beta=X_c$.  This
process in not unique --- each bead stabilizes to exactly
$d= A_i\cdot X_c$ stable maps.

\medskip

The following definition will help make precise which maps occur in a
necklace.
\begin{defn} {\em Let  $\M^\nu_{d,A}$ denote the moduli space of perturbed maps
$\phi:(S^2,x,y,u)\ra S^2\ti X$ representing $A$ in homology and
such that $\phi$ has a contact of order $d$ to $X_a$ at $x$, a contact of
order $d$ to $X_b$ at $y$, and $\phi(u)\in X_c$
(where $d=A\cdot [X_a]$).}
\end{defn}

The moduli space $\M^\nu_{d,A}$ comes with an evaluation map
$$
\ev:\M^\nu_{d,A} \ra X\ti X
$$
by
sending $(\phi,x,y)$ to $(p_X\phi(x),p_X\phi(y))$.  This extends to  an
evaluation map
\bear\label{evalmap}
\ev:\M^\nu_{d,A_i}\ti\dots\ti\M^\nu_{d,A_k} \ra (X\ti X)^{k}.
\eear
The condition that each map $(\phi_i,x_i,y_i)$ in the necklace connects to
the next is
$$
f(p_X\phi_i(x_i))=p_X\phi_{i+1}(y_{i+1})
$$
(and similarly the last map connects to the first).  These conditions say
that  $\ev((\phi_i,x_i,y_i)_{i=1}^k)$ lies on
the ``cyclic graph'' of $f$:
\bear
G_f=
 \left\{ \;(x_1,y_1,\dots,x_{k},y_{k})\;|\;
f(x_i)=y_{i+1}\mbox{ for $i=1,\dots,k-1$ and }f(x_k)=y_1\ \right\}.
\label{cyclicgraphdef}
\eear
Thus  the moduli space of stabilized necklaces is  parameterized by
\bear
\ev^{-1}(G_f)=\bigsqcup\;\ev (\M^\nu_{d,A_1}\ti\dots\ti\M^\nu_{d,A_k})\cap G_f
\label{firstevGf}
\eear
where  the union is over  collections $A_1,\dots
A_k$ of elements $A_i=ds+\alpha_i$ in $H_2(S^2\ti X)$ satisfying
(\ref{5.sumalpha=0}).  Note
that the  evaluation map is
orientation preserving on moduli spaces of spheres so  pushing forward by the
evaluation map (\ref{evalmap})  is the same as evaluating by pullback.

\medskip

Thus far we have assumed that the bubble tree limit in $X^0_f$ has no
components entirely in the singular
set.  To achieve such a condition in general it is necessary to renormalize
in a slightly different way.
Consider the original sequence of maps $f_{\ep}$ into $X^{\ep}_f$.  Instead
of  simply letting $\ep\to 0$, we
choose  $r$ and sequences $0<\ep^1_i<\ep^2_i<\dots <\ep^r_i$ all
approaching 0 as $i\to\infty$,
renormalize the annuli
$A^\ell_a=N_a(\ep^{\ell+1})\setminus N_a(\ep^{\ell})$, and collapse their
boundaries by a ``symplectic cut'' (cf. [IP2]).  The resulting bubble tree
limit goes not into $X^0_f$,  but instead into the singular space
\bear
\label{X0r}
X^{r}_f= (S^2\ti X)_0\cup (S^2\ti X)_1\cup \dots \cup (S^2\ti X)_r\  /\sim
\eear
where $\sim$ identifies $X_b$ on one copy of $S^2\ti X$ to $X_a$ on the
next by the identity map, and identifies $X_b$ on last copy to $X_a$ on the
first by the map $f$.  This renormalization insures that the limit
\bear
\begin{array}{l}
\mbox{(a)\  has no components entirely contained
in the singular set of $X_f^{r}$, and}\\
\mbox{(b)\  the total homology class in
each copy of $S^2\ti X$ besides the first one}\\
\mbox{\qquad {\em is not  a multiple of $s$ }
(otherwise there was no need to renormalize)}.
\end{array}
\label{rcondtions}
\eear
We can then remove double points,   lift, and remove  singularities  to
obtain a map
$\bigsqcup \ov B_i\to S^2\ti X$ exactly as in (\ref{5.liftedBT}). The
rest of the discussion carries through.  The necklaces then lie in
$\ev(G^r_f)$, where $G^r_f$ is the cyclic
graph of the $r+1$ string $(id,id,\dots id, f)$, that is, the subset of
$(X\ti X)^{(r+1)k}$ defined as in
(\ref{cyclicgraphdef}) where the identifications are in the order
$(id,\dots,id,f;\ \dots \   ;id,
\dots id,f)$.

Thus  the moduli space of stabilized necklaces is  parameterized by the set
analogous to (\ref{firstevGf})
with $G_f$ replaced by $G_f^r$ for $r=0,1,\dots$. Because each unstable
bead stabilizes in $d$ ways,
the ``unstable'' necklaces of Figure 4.2  are parameterized by the set
$$
{1\over d^{k-1}}\; \bigsqcup_{r,\cal A}\;
\ev (\M^\nu_{d,A_1}\ti\dots\ti\M^\nu_{d,A_k})\cap G_f^r.
$$
where  the union is over all $r$ and over the set ${\cal A}$ of all
collections $A_1,\dots
A_k$ of elements $A_i=ds+\alpha_i$ in $H_2(S^2\ti X)$ satisfying
(\ref{5.sumalpha=0}) and (\ref{rcondtions}b).

\bigskip

Going the other way, a gluing argument shows that each necklace in
$\ev^{-1}(G_f^r)$ can be perturbed to exactly $d^k$ approximate
$(J,\nu)$-holomorphic tori in $X_f^\ep$, and that each of
these can be uniquely corrected to a perturbed holomorphic map.  The
approximate maps are obtained by an essentially canonical ``rounding off''
procedure applied to each of the $k$ points where the beads
intersect.  In general, whenever two stable maps have a contact of order
$d$, we obtain $d$ approximate maps; these have the same image, but have
distinct complex structures on their domains.  The details of this
construction can be found in [IP2].

\medskip

Altogether, each element of $\ev^{-1}(G_f^r)$ corresponds to
$d^k/d^{k-1}=d$ distinct $(J,\nu)$-holomorphic tori.
This,  combined with the main result of [IP2]
therefore implies the following.
\begin{proposition}
\label{RT convergence}
For generic $(J',\nu')$ on $X_f$ and generic $(J,\nu)$ on $S^2\ti X$,
there is an oriented cobordism between the moduli
space (\ref{5.3}) defining $RT_m$ for $(J',\nu')$  and a
covering of the space of $(J,\nu)$-necklaces:
\best
\M_{mT}(T^2,X_f)\ \sim \  \ma\bigsqcup_{r,\;{\cal A}}
\; \ma\bigsqcup_{i=1}^d\;
\ev (\M^\nu_{d,A_1}\ti\dots\ti\M^\nu_{d,A_k})  \cap G_f^r
\eest
where the last union is over $d$ identical disjoint copies.
\end{proposition}

Thus the Gromov invariant we seek can be  written as a sum of
homological intersections:
\bear
RT_m\ =\ \sum_{r,\;{\cal A}}\ d\cdot
\left[\ev (\M^\nu_{d,A_1}\ti\dots\ti\M^\nu_{d,A_k})\cap G_f^r\right].
\label{5.intersectioninSX}
\eear
This formula expresses
$RT_m$ (an invariant of $X_f$) as a count of curves in  $S^2\ti X$.   By
construction, the righthand side counts perturbed holomorphic curves for a
rather special set of pairs
$(J,\nu)$, namely the ones for which the restrictions of $J$ and
$\nu$ to $X_a$ and $X_b$ agree under $f:X_a\to X_b$.  But in fact the
righthand side is a Ruan-Tian invariant on
$S^2\ti X$, and hence  can be computed using {\it any} generic $(J,\nu)$ on
$S^2\ti X$.  We will exploit this
by fixing a product almost complex structure $J_0$ on $S^2\ti X$,  showing
that the pair
$(J_0,0)$ is generic, and then explicitly computing the invariant for
$(J_0,0)$.

\medskip

The first step is to show that the only necklaces which contribute to the
invariant are those where every bead represents a multiple of
$s$ in homology.  The proof hinges on the observation that the symplectic form
 $\w$ on $X$ pullsback to a 2-form $\bar{\w}$ on $S^2\ti X$; this
satisfies $\bar{\w}|_{X_a}=f^*\bar{\w}|_{X_b}$, so descends to a form $\w_f$
on $X_f^r$.  These forms
determine cohomology classes $[\w]$, $[\bar{\w}]$ and $[\w_f]$.

\begin{lemma}
\label{J_0isOK}
For all $(J,\nu)$ close to $(J_0,0)$, the only non-vanishing intersections
in (\ref{5.intersectioninSX}) are
those with $r=0$ and  $A_i=d\cdot s$ for all $i$. Thus
\bear
\label{RTmFormula}
RT_m\ =\ \sum_{m=dk}d\cdot
\left[\ev (\M^\nu_{d,ds}\ti\dots\ti\M^\nu_{d,ds})\cap G_f\right].
\eear
\end{lemma}
\pf  Suppose that for some choice $(A_1,\dots A_k)\in {\cal A}$ the
intersection (\ref{5.intersectioninSX}) is
non-empty for all $(J,\nu)$ near $(J_0,0)$.  Choose a sequence
$(J_j,\nu_j)\to (J_0,0)$ and, for each
$j$,  a bubble tree $B_j=\cup B_{ij}$ in the intersection
(\ref{5.intersectioninSX}) with the component
$B_{ij}$ representing $A_i$ in homology for each $j$.   The Bubble Tree
Convergence Theorem applies to each
component, showing that a subsequence of  $\{B_{ij}\}$ converges to a bubble
tree limit $B_{i\infty}$ in $X_f^r$ and lifts to a bubble tree
$\bar{B}_{i\infty}$ in $S^2\ti X$.  Let
$\{\bar{C}_{i,\ell}\}$ be the components of this limit and write the
homology class $[\bar{C}_{i\ell}]$ as $d_{i\ell}s + \alpha'_{i\ell}$ as we did
before  (\ref{5.sumalpha=0}).  Then
for each $i$ we have $A_i=\sum_{\ell}[\bar{C}_{i\ell}]=ds+\alpha'_i$.
Furthermore, each $\bar{C}_{i\ell}$ is
$(J_0,0)$-holomorphic for the product structure $J_0$, so the union
$$
{\cal C}=\bigcup_{i,\ell}\bar{C}_{i,\ell}
$$
projects to a holomorphic curve $p({\cal C})$ in $X$ representing $\sum_{i,\ell}
\alpha'_{i,\ell}=\sum_i\alpha'_i$ with
$$
\sum\mbox{Area}(p(\bar{C}_{i,\ell}))   \ =\ \mbox{Area}(p({\cal
C}))\ =\  [\w]\cdot [{\cal C}]\ =\
[\bar{\w}]\cdot[{\cal C}].
$$
On the other hand, $\bar{\w}$ descends to $\w_f$ under the projection $\pi$
from $S^2\ti X$ to
$X_f^r$, and   by (\ref{5.sumalpha=0}) $\pi_*[{\cal C}]$
represents a multiple of $T_0$ in $H_2(X_f^r)$.  Thus
$$
[\bar{\w}]\cdot[{\cal C}]\ =\ [\w_f]\cdot \pi_*[{\cal C}]=0.
$$
We conclude that
$p({\cal C})$ is trivial, so all the
$\alpha'_{i,\ell}$ are zero.  Therefore $A_i=ds$ for each $i$, and hence
$r=0$ by  property (\ref{rcondtions}b).
\qed

\medskip

\begin{lemma}
\label{productgeneric}
The choice $(J,\nu)=(J_0,0)$ is generic
for the class $ms$ and for this choice there is an oriented identification
$$
\ev (\M^0_{d, ds}\ti\dots \ti\M^0_{d, ds})=\De^k
$$
where $\De^k=\De\ti\dots\ti\De$ is the multidiagonal in $ (X\ti X)^{k}$.
\label{lemmacount1}
\end{lemma}
\pf  Modulo diffeomorphisms, we may assume that the marked points on $S^2$ are
$a=x=0,\;c=u= 1$ and $b=y=\infty$. Then, for the choice $(J_0,0)$,
\best
\M^0_{d, ds}=\{ \phi:S^2\ra S^2\ti X\;|\; \phi(z)=(z^d,p),
\mbox{ for some fixed } p\in X\;\},
\eest
so $\M^0_{d,ds}$ is canonically identified with $X$. The linearization of the
$J$-holomorphic map equation at $\phi\in\M^0_{d,ds}$ is then the
$\overline{\partial}$ operator on the pullback normal bundle to the
image. Then $(J_0,0)$ is generic because at  $\phi\in\M^0_{d,ds}$ the
 cokernel of the linearization vanishes:
\best
\mbox{ cok }\ov{\partial}=H^{1}(S^2,\phi^*T X)=
H^{0,1}(\P^1,{\cal O})\otimes T_pX=0.
\eest
The moduli space is then oriented by the complex structure on
\best
T_{\phi}\M^0_{d,ds}\ =\ \mbox{ker } \overline{\partial}\ =\
H^0(S^2;\phi^*TX)\ =H^{0}(\P^1,{\cal O})\otimes  T_pX= \ T_pX.
\eest
Thus the identification $\M^0_{d,ds}\cong X$ respects orientation.  Under the
corresponding identification of $\M^0_{d,ds}\ti\dots\M^0_{d,ds}$ with
$X\ti\dots\ti X$, the evaluation map (\ref{evalmap}) becomes  the map
\best
\ev(p_1,\dots,p_{k})=(p_1,p_1,\dots,p_{k},p_{k} )
\eest
whose image is $\Delta^k$.
\qed

\medskip

By Lemma \ref{lemmacount1}, the intersection in  (\ref{RTmFormula})
for the product complex structure becomes
\bear
\label{5.Grintersection}
\ev (\M^0_{d,ds}\ti\dots \ti\M^0_{d,ds})\cap G_f=\De^k\cap G_f
\eear
for each $d$ and $k$.  This is directly related to the Lefschetz number of $f$.

\begin{lemma}\label{5.alglemma} The homological intersection
$\De^k\cap G_f\ $ is $\ L(f^k)$.
\end{lemma}
\pf   After  perturbing  as in Lemma \ref{transversality},  $L(f^k)$ is,
by definition, the sum of the fixed points of $f^k$, each oriented by the
sign of
$\mbox{det }(df^k-I)$.  It is easy to see that there is a one-to-one
correspondence between
the fixed points of $f^m$ and the points of $\Delta^k\cap G_f$.  To verify
that the signs agree,
consider the isomorphism $B$ of $TX^k\oplus TX^k$
given by the $2k\ti 2k$ block matrix
$$
B=\l(\begin{array}{cc}I&I\\I&A\end{array}\r)
\qquad
\mbox{where}
\qquad
A=
\l(\begin{array}{ccccc}
0& df&0  && 0\\
0&0&df&\dots&0\\
&\vdots& &  &df\\
df&0&0&&0\\
\end{array}\r).
$$
Then $B$ takes the first factor $TX^k$ into $T\Delta^k$ and the second
factor to $TG_f$.
The sign of an intersection point in $\Delta^k\cap G_f$ is therefore the
sign of
$$
\mbox{det }B=\mbox{det }(A-I)=\mbox{det }(df^k-I)
$$
(to obtain the last equality  expand along the top row).
\qed

\bigskip

Combining Lemma \ref{5.alglemma} with equations (\ref{RTmFormula}) and
(\ref{5.Grintersection}) gives
\bear
RT_{m}&=&\ma\sum_{m=kd}d\; L(f^{k}).
\label{RT=Lf}
\eear
 Theorem \ref{theorem5.1} now follows by simple algebra.  From equations
(\ref{5.2}) and
(\ref{RT=Lf})
\best
\log Gr^T(t)\ =\ \ma\sum_{m=1}^\infty{RT_m\over m} \log  F(t^m)=
\ma\sum_{m=kd} {d L(f^{k})\over kd} \log  F(t^{kd})=
\ma\sum_{k=1}^\infty {L(f^k)\over k}\ma\sum_{d=1}^\infty \log F(t^{kd}),
\eest
Using  relation (\ref{prodF=t}) this simplifies to
\best
\log Gr^T(t)\ =\ \ma\sum_{k= 1}^\infty {L(f^k)\over k} t^k,
\eest
which exponentiates to give $Gr^T(t)=\zeta_f(t)$.  This completes the proof
of Theorem \ref{theorem5.1}.

\bigskip

\setcounter{equation}{0}
\section{\bf Gluing in $X_f$}\label{Gluing in}
\bigskip

Thus far,  we have constructed the symplectic mapping tori $X_f$ and
computed their Gromov invariants.  In this section we will   view the
 $X_f$ as {\it operations} on the set of
symplectic four-manifolds.  While the basic construction works in all
dimensions, it is particularly
interesting in dimension four, where the mapping tori are related to knots,
and  the Gromov invariants
are related to the Seiberg-Witten invariants.

\medskip

Taking  $X$ a genus $g$ Riemann surface and $f$ a symplectomorphism of $X$
such that
\bear\label{det=pm1}
det(I-f_{*1})=\pm 1
\eear
gives symplectic 4-manifolds $X_f$ which is homology $S^2\ti T^2$ with a
well-defined section class
$T\subset X_f$  represented by a symplectically embedded torus, and a
fiber class $F$ represented by a symplectically embedded genus $g$ surface.
Both $T$ and $F$ have square zero, and $F\cdot T=0$.
\begin{theorem}For  $f$  as above, the  Gromov series of $X_f$  is
\bear
Gr(X_f)= \frac{\mbox{det}(I-t_T\cdot f_{*1})}{(1-t_T)^2}.
\label{eq6.1}
\eear
\label{thm6.1}
\end{theorem}

\pf Theorem \ref{introThm1} shows that the {\it partial} Gromov series
$Gr^T(X_f)$ has this form.  Thus we
need only show that $Gr^T(X_f)=Gr(X_f)$, that is, that any class $A\in
H_2(X_f)$  with $RT(A)=0$ is a
multiple of $T$.  Write such a class as $A=\alpha F +\beta T$.  By choosing
a sequence of generic pairs
$(J_i,\nu_i)$ converging to $(J,0)$ and applying the Bubble Tree
Convergence Theorem, we see that $A$ can be
represented by a $J$-holomorphic curve, so $\alpha = A\cdot T\geq 0$.
Then, by (\ref{1.defdegzero}) and
Lemma \ref{canonical class}  below, we have $g-1=\kappa\cdot A=-2\alpha$.
Thus $\alpha=0$ and $A$ is a multiple of $T$.
\qed

\medskip

Now suppose  $(Z,\w)$ is a
symplectic 4-manifold and  $F\subset Z$ is  a symplectically embedded
torus with trivial normal bundle (``symplectic torus fiber'').
Then for each symplectomorphism $f$ of the Riemann surface $X$ we can form
the symplectic fiber sum
\bear\label{Zf}
Z(f)\ =\ Z\ma\#_{F=T} X_f
\eear
by removing tubular neighborhoods of $T\subset X_f$ and $F\subset Z$ and
gluing. Thus each  symplectomorphism $f$ gives an operation $Z\mapsto Z(f)$
on the set of symplectic 4-manifolds with symplectic torus fibers. In
particular, for each fibered knot $K$ we get an operation $Z\mapsto Z(K)$ by
taking $f$  to be the monodromy of $K$.

\begin{theorem} Let $(Z,\w)$ be a symplectic 4-manifold that contains a
symplectically torus fiber $F$. If $f$ satisfies (\ref{det=pm1}) then the
 partial Gromov invariant for the modified manifold (\ref{Zf}) is
\bear\label{grZf}
Gr^F_{Z(f)}\ =\ Gr^F_Z\cdot \mbox{det}(I-t f_{*1}).
\eear
If (a) the linking  circle of $F$ is homotopic to zero in
$Z\setminus F$ and (b) $f$ is the monodromy of a fibered knot $K$ in a
 homotopy 3-sphere, with Alexander polynomial $A_K$ then $Z(K)$ is homotopic
to $Z$
\bear\label{ZK=Z}
Z(K)\sim Z
\eear
via a homotopy that preserves $F$ and
\bear\label{grZK}
Gr^F_{Z(K)}\ =\ Gr^F_Z\cdot A_K(t_F).
\label{6.grZ'}
\eear
\label{glueXftoY}
\end{theorem}
\pf Combining (\ref{1.gluinglemma}) with (\ref{eq6.1}) gives (\ref{grZf}) and
(\ref{grZK}). If $X_f$ arises from a fibered knot then $X_f=M_f\ti S^1$,
where by (\ref{Mf-tau})
$\pi_1(M_f\setminus \{\tau\})=\pi_1(S^3\setminus K)=\Z$ is
generated the linking circle of $\tau$. The homotopy equivalence
(\ref{ZK=Z}) then follows from the Van Kampen Theorem, applied to the
 decompositions $(Z\setminus F)\cup (X_f\setminus T)$ and
$(Z\setminus F)\cup (S^2\ti T\setminus T)$. \qed

\medskip

Note how knots enter this picture.  For {\em any} symplectomorphism $f$ of
Riemann surfaces that satisfies (\ref{det=pm1}), we can form $X_f$
and modify $Z$ to $Z(f)$, and get formula (\ref{grZf}) with
$A(t)=\mbox{det}(I-tf_{*1})$.  But $Z$ and
$Z(f)$ will have different homotopy types  unless $\pi_1(M_f)$ is normally
generated by $\tau$, and it is this condition that leads us to knots.
\medskip

>From (\ref{grZK}) and Lemma \ref{l2.5} we obtain:
\begin{cor}Let $Z$ be a symplectic 4-manifold  with a symplectic torus fiber
 $F$ such that the linking circle of $F$ is homotopic to zero in
$Z\setminus F$. Then $Z$ is homotopic to infinitely many distinct symplectic
4-manifolds $Z(K)$.
\end{cor}
\bigskip

In the remainder of this section we will apply Theorem \ref{glueXftoY} to
elliptic surfaces.  Thus we take
$Z=E(n)$ and consider the manifolds
$$
E(n;K)=E(n)\#X_K.
$$
For most of these, we can compute the full Gromov series. The following
result appeared as Theorem
\ref{introThm2} of the introduction.
\begin{theorem}
\label{fullGr=Grt}
For $n\neq 1$, the  Gromov series of $E(n,K)$ is
\begin{equation}
\label{GrofEnK}
Gr(E(n,K))\ =\ (1-t)^{n-2}\ A_K(t)
\end{equation}
where the variable $t$ corresponds to the fiber class, and $A_K$ is the
Alexander polynomial of the knot. For $n>1$ the same formula gives the
Seiberg-Witten series:
\begin{equation}
\label{SWofEnK}
SW(E(n,K))\ =\ (1-t)^{n-2}\ A_K(t).
\end{equation}
\end{theorem}
\pf When $n=0$ (\ref{GrofEnK}) is the same as (\ref{eq6.1}), so take $n>1$.
Let $C\in H_2(E(n,K))$ be a class
with $RT(C)\neq 0$.  We will show that $C$ is a multiple of the fiber class
$F$, and then (\ref{GrofEnK}) will follow from
Theorem \ref{glueXftoY} and (\ref{1.GTEn}).

    Since  $n>1$ we can write $E(n,K)$ as a symplectic fiber sum
$$
E(n,K)= E(1) \#_{F_1} E(2) \#_{F_2} \cdots  \#_{F_{k-1}} E(2) \#_{F_k} X_K
$$
where there is at least one copy of $E(2)$ and at most one copy of $E(1)$.
This fiber sum depends on
parameters $\ep_i$,  the radii of the tubular neighborhoods of the fibers
$F_i$ removed in its construction.
Choose a sequence $\ep_i\to 0$ and a sequence of pairs $(J_i,\nu_i)$
generic on $E(n,K)_{\ep}$ with $\nu_i\to 0$
and $J_i$ converging to  a generic K\"{a}hler structure  $J_0$ on each
$E(2)$.  Then fix a sequence $C_i$ of
$(J_i,\nu_i)$-holomorphic curves representing $C$.  In the limit, these
converge to a bubble tree curve $C_0$ in
the singular space obtained by identifying $E(1)$, $E(2)$, and $X_K$ along
the fibers $F_i$ as above.

   Since a generic K\"{a}hler structure  on $E(2)$ admits no holomorphic
curves, the  curve $C_0$ is a disjoint
union of curves $A\subset E(1)$ and $B\subset X_K$ which, moreover, do not
intersect the fiber.
  Since the canonical divisor is a multiple of the fiber in
both $E(1)$ and $X_K$ (Lemma (\ref{canonical class}) below), the adjunction
inequality implies that $A\cdot A\geq 0$ unless at least one component of
$A$ is an exceptional sphere, and that
cannot be true since in $E(1)$ all exceptional spheres intersect the fiber.
We similarly conclude that $B\cdot B=0$.

We now have classes $A,F$ in $E(1)$ with $A\cdot A\geq 0$ and  $F\cdot
F=A\cdot F=0$.  Since $E(1)$ has $b^+=1$,
the ``Lightcone Lemma'' ([LL]), implies that $A$ is a multiple of $F$.
Since $b^+(X_K)=1$ the same argument
applies to $B$, finishing the proof of (\ref{GrofEnK}).

To connect this with the Seiberg-Witten series, note that the manifolds
$E(n,K)$ have $b^+>1$ when
$n>1$.  For such manifolds, Taubes proved that the Seiberg-Witten invariants
are equal
to the
Taubes-Gromov invariants [T3] and that both of these invariants are
nontrivial only for zero-dimensional moduli spaces  (i.e.
unconstrained curves) [T2].  On the other hand, we have just shown that the only
nontrivial unconstrained Ruan-Tian invariants for $E(n,K)$ correspond to
tori representing
multiples of $T$. But for these, the Ruan-Tian series defined above agrees
with the one defined in [IP1], equation (3.3), which in turn is equal to the
Taubes-Gromov series by the main result in [IP1].
\qed

\bigskip

Because the Seiberg-Witten invariants depend only on the differentiable
structure, Theorem \ref{GrofEnK}
implies that  $E(n,K)$ and $E(n,K')$ are distinct differentiable manifolds
whenever $n>1$ and $K$ and $K'$
have different Alexander polynomials.  Thus by Burde's Theorem there are
infinitely
many distinct symplectic 4-manifolds homeomorphic to $E(n)$.

\medskip

\noindent{\bf Remark}  For the exceptional case $E(1,K)$, we still have
that the {\it partial} Gromov series
$Gr^T$ is given by $(1-t)^{n-2}\ A_K(t)$, but there are other classes ---
such as the exceptional curves ---
which contribute to the full Gromov invariant.  On the other hand,
Fintushel and Stern showed that the
Seiberg-Witten series of $E(1,K)$ is $(1-t)^{n-2}\ A_K(t)$.  Thus in this
case the Gromov series defined
in section 1 carries more information than the SW series.

\bigskip

The above proofs required explicit expressions for the canonical
divisors.   These are
calculated in next lemma.

Let $F$ denote  the `fiber class' of $E(n;K)=E(n)\#X_K$, which can be
thought of as either the fiber class of $E(n)$ or the section class $T$
of $X_K$. Then
\begin{lemma}
The canonical divisor of $X_f$ is $\kappa=(2g-2)T$ and the canonical
divisor of $E(n,K)$ is $\kappa=(2g-2+n)F$
\label{canonical class}
\end{lemma}
\pf For both parts we will use the adjunction formula:  if $A\in H_2(X)$
is represented by a symplectically embedded curve of genus $g$ then
\bear\label{adj}
\kappa\cdot A = 2(g-1)-A\cdot A.
\eear
Recall that  $H_2(X_f)$ is generated by classes $F$ and  $T$, both
represented by a symplectically embedded
curves of square zero and genus $g$ and 1 respectively, such that
$F\cdot T=1$.  Then (\ref{adj}) uniquely determines $\kappa=(2g-2)T$.

Next, fix a fiber $F_0$ and a holomorphic section $\sigma$ in $E(n)$;
$\sigma$ is symplectic with genus $0$,
$\sigma\cdot \sigma=-n$, and $\sigma\cdot F_0=1$.  Similarly, fix a section
$T$ of $X_K\to T^2$.  For each
$p\in T$ there is a fiber $X_p$ (a symplectically embedded copy of $X$)
with $X_p\cap T=\{p\}$ and
with self-intersection zero. We can then identify $F_0$ with $T$ and form
the fiber sum $E(n,K)$, simultaneously gluing
$\sigma$ to the appropriate $X_p$, obtaining a  symplectically
embedded genus $g$ curve $s$ in $E(n,K)$ (cf. [G] Theorem 1.4) with
$s\cdot s=-n$.

Now  let $\kappa$ be the canonical divisor of $E(n,K)$.   Then $s\cdot s=-n$,
$s\cdot F=1$, and $F\cdot F=0$, so we can complete $\{s,F\}$ to a basis
$\{s,F, a_i\}$ of $H_2(E(n,K))$ with $a_i\cdot
F=0$ for all $i$.  But any class $a$ with $a\cdot F=0$ can be represented
by the disjoint union of cycles
$C\subset E(n)$ and $C'\subset X_K$ with $C\cdot F=0$ and $C'\cdot t=0$.
Along $C$, the canonical bundle of
$E(n,K)$ is the same as the canonical bundle of
$E(n)$, so $\kappa\cdot C=(n-2)F\cdot C=0$.    Similarly, $\kappa\cdot
C'=(2g-2)T\cdot C=0$.  Thus $F$ and all
the classes $a_i$ are perpendicular to $\kappa$. This and the equation
$\kappa\cdot s = 2g-2+n$ (from the
adjunction formula) uniquely determine $\kappa$ to be $(2g-2+n)F$.
\qed

\medskip

Finally, the following result verifies the assertion made before Theorem
\ref{introThm2}.

\begin{prop}
For any fibered knots $K,\; K'$ there is a homeomorphism between $E(n,K)$
and $E(n,K')$
preserving the fiber class, the periods of $\w$, and the Stiefel-Whitney
class $w_2$. When  $K$  and $K'$
have the same genus this homeomorphism also preserves the canonical class.
\label{6.last}
\end{prop}

\pf  Since the elliptic surfaces $E(n)\to {\Bbb P}^1$ admit sections, the
linking circle of a  fiber $F$ is
contractible in $E(n)\setminus F$.  Hence by Theorem \ref{glueXftoY} and
the classification of
simply-connected topological 4-manifolds there is a homeomorphism between
$E(n,K')$ and $E(n,K)$ that preserves $F$.  This homeomorphism preserves the
canonical divisor $\kappa=(2g-2+n)F$ when
$K$ and $K'$ have the same genus, and preserves $w_2$(the mod 2 reduction
of $\kappa$) regardless of the
genus.  It remains to show that we can choose the symplectic structures so that
the periods of $\w$ match under this homeomorphism.

Recall that $X_K$ is the symplectic mapping cylinder (\ref{eq00.1}) of the
monodromy maps $f$ of $K$.  Then
$H_2(X_K)$ is generated by the section class $T$ with $\w(T)=1$ and the
dual class $F$ with
$\w(F)=\mbox{vol}(X)$ where $X$ is the fiber of the knot.  Thus the periods
of $\w$ on $X_K$ and
$X_{K'}$ match.

Now let $\{s,F, a_i\}$ be the basis of $H_2(E(n,K))$ used in Lemma
\ref{canonical class}.  Since the
classes $F$ and $a_i$ have representatives that are disjoint from the
gluing region (a neighborhood of one
fiber), they are preserved under the homeomorphism and have the same
periods.   The period of the remaining
generator $s$ also matches since it is formed in $E(n,K)$ and $E(n,K')$ by
applying the same construction to
curves of equal volume.
\qed

\bigskip
\setcounter{equation}{0}
\section{\bf Symplectic Structures in Dimension Six}
\bigskip

It is difficult to distinguish differentiable structures on 4-manifolds,
but in higher dimensions
differentiable structures are well-understood via surgery theory.  The
simplest way to move to higher dimensions
is to take the product with a 2-sphere.  In this context, one can show that
if $X$ and $Y$ are
smooth oriented 1-connected  4-manifolds which are homeomorphic and have
the same characteristic classes, the  6-manifolds $X\ti S^2$ and $Y\ti S^2$
are diffeomorphic.  We
will use that fact to construct
families of manifolds that admit infinitely many deformation classes of
symplectic structures.  Previously, Ruan
and Tian [RT] have found such manifolds.

 The examples described here are remarkable for their simplicity.
Furthermore, the differences in the
symplectic structures cannot be detected by {\it any} of the    ``classical
invariants'' mentioned in the introduction: the  diffeomorphism type,
the  cohomology class of the symplectic form, and the
homotopy class $[J]$.

\begin{theorem}
The differentiable 6-manifolds $K3\ti S^2$ and, more generally,
$E(n)\ti S^2$ have infinitely many deformation classes
of symplectic structures, all with the same $[J]$ and $[\w]\in
H^2_{DR}$.
\label{7.thm1}
\end{theorem}
We will prove this after giving three lemmas. The proof actually
shows more: for each integer $g> 1$ there are infinitely
many deformation classes of symplectic structures whose canonical
divisor is $\kappa=\kappa_{E(n)\ti S^2} -2g F\ti [S^2]$.

\medskip

Theorem \ref{7.thm1} is a result of comparing the symplectic manifolds
$E(n,K)\ti S^2$ for different knots
$K$.  The first lemma shows that these cannot be distinguished by  the classical
invariants. Let  $F$ denote  the fiber class in
$E(n,K)$ and let $F'=F\ti\{\mbox{pt}\}$ be the corresponding class in
$E(n,K)\ti S^2 $.

\begin{lemma}

For any two fibered knots $K,K'$, there is a diffeomorphism $E(n,K)\ti S^2 \to
E(n,K')\ti S^2$ that preserves $[\w]$.  When  $K$  and $K'$ have the same genus
this diffeomorphism also preserves $[J]$.
\label{lemma7.1}
\end{lemma}

\pf By Proposition \ref{6.last} there is a homeomorphism $\phi: E(n,K) \to
E(n,K')$ which preserves
$[\w]$  and $w_2$ (and also $\kappa$ when $g$=genus$(K)$ equals
$g'$=genus$(K')$). It also preserves $p_1$
because for 4-manifolds $p_1$ is three times the signature and because
signatures add under fiber sums, so
$\sigma(E(n,K))=\sigma(E(n)) +\sigma(X_K)=\sigma(E(n))$.

Now  $E(n,K)\ti S^2$ and $E(n,K')\ti S^2$ are homeomorphic by
$\psi=\phi\times id.$   Since the symplectic form, the Pontryagin class,
and the canonical
class  all respect the product structure (e.g.
$\kappa=\pi_1^*\kappa_{E(n,K)} +\pi_2^*\kappa_{S^2}$), $\psi$ preserves
$[\w]$ and $p_1$, and
preserves $\kappa$ when $g=g'$.  We can then apply the result mentioned
above:  if
$X,Y$ are  closed,
1-connected, oriented  smooth 6-manifolds with torsion-free homology, and
$X\to Y$ is a homeomorphism
that  preserves orientation, $w_2$ and $p_1$, then there is a
diffeomorphism $X\to Y$ that induces the
same map on homology (see [J]).  Thus we get the desired diffeomorphism,
and it preserves $[\w]$ and preserves $\kappa$ when $g=g'$.  But for
6-manifolds the homotopy class of the
tamed almost complex structures is determined by
$\kappa$ ([Ru]), so $[J]$ is also preserved when $g=g'$.
\qed

\medskip

  Lemma \ref{lemma7.1} shows that we can
regard  $E(n,K)\ti S^2$ as a symplectic structure on $E(n)\ti S^2$, and
that two of these structures have
the same classical invariants whenever the corresponding knots have the
same genus.   We
will complete the proof of Theorem
\ref{7.thm1} by showing that these symplectic structures can be
distinguished by Gromov invariants.  Recall
from section 1 that when $n=3$ there is a separate Gromov series $Gr_g$ for
each genus $g$.  We will compute
the genus 1 Gromov series and show that
$$
Gr_1(E(n,K)\ti S^2)\neq Gr_1(E(n,K')\ti S^2).
$$
The first observation is that the  genus 1 Gromov series  are well-behaved
with respect to products.

\medskip

\begin{lemma}\label{GrX*S2} Let $X$ and $Y$ be symplectic manifolds and
consider  $Z=X\times Y$ with the
product symplectic structure.  Let $\chi(Y)$ be the euler characteristic of
$Y$.  Then
$$
Gr_1^{A\ti \{pt\}} (X\ti Y) =\left[ Gr_1^A(X)\right]^{\chi(Y)}.
$$
\end{lemma}
\pf  By the definition of the genus 1 Gromov series, this follows from the
equality
\bear
RT_{A\ti \{pt\},1,X\ti Y}(\al_1\ti [Y],\dots,\al_k\ti [Y])=
\chi(Y)\cdot RT_{A,1,X}(\al_1,\dots,\al_k).
\label{rtX*S2}
\eear
for any $A\in H_2(X,\Z)$ and constraints $\al_1,\dots, \al_k\in H_*(X,\Z)$.
Let $\M$ denote the moduli space
corresponding to the lefthand side of (\ref{rtX*S2}), that is, the set of
genus one perturbed holomorphic
maps from a torus to $Z$ representing $A\ti \{pt\}$ and passing through the
given constraints.  Similarly, the
$\M'$ denote the moduli space corresponding to the righthand side of
(\ref{rtX*S2}).  Choose a generic
$(J_0,\nu_0)$ on $X$, a generic $J_1$ on $Y$, and consider
$(J_0\oplus J_1,\nu_0)$ on $Z$; with this choice $\nu=0$ in the tangent
direction to $Y$.

The moduli space of $(J_0\oplus J_1,\nu_0)$-holomorphic maps in $Z$  consists of
maps $f:T^2\ra Z$ whose projection $\pi_X f$ on the $X$ factor
is an element of $\M'$ while $\pi_Y f$ is a constant map.  Thus $\M$ is
canonically diffeomorphic to a product
of copies  of $Y$, one for each map in  $\M$ and each with a sign. More
precisely,
\bear
\M\,\cong\, \underbrace{Y\sqcup\dots\sqcup Y}_p\,\sqcup\,
 \underbrace{-Y\sqcup\dots\sqcup -Y}_n \;,
\label{7.copiesofY}
\eear
with  $RT_{A,1,X}=p-n$.

The choice $(J_0\oplus J_1,\nu)$ is not generic, so the invariant must be
calculated as the Euler class of
the Taubes obstruction bundle over $\M$.   The fiber of the obstruction
over $f\in\M$ is the cokernel of
 the linearization  $D_f$ of the perturbed holomorphic equation.
In our case, because of the  product structure and because $\pi_Y f$ is the
constant map to some point $p\in
Y$,
\best
\mbox{coker }D_f=\mbox{coker }D_{\pi_1 f}\oplus \mbox{coker }D_{\pi_2 f}=
0\oplus H^{0,1}(T^2,(\pi_Y f)^*(TY))\cong T_{p} Y
\eest
($D_{\pi_1 f}$ has no cokernel since $(J_0,\nu_0)$ is generic on $X$).
Thus the obstruction bundle is
canonically isomorphic to the tangent bundle of $Y$  over each copy of $Y$
in (\ref{7.copiesofY}).
Consequently, when we perturb to a generic $(J,\nu)$, each copy of $Y$ in
(\ref{rtX*S2}) gives rise to
exactly $\chi(Y)$  perturbed holomorphic maps with sign, for a total of
$\chi(Y)\cdot RT_{A,1,X}$.
\qed

\medskip

\begin{lemma}
If $X=E(n,K)$ with fiber class $F$, then
\begin{equation}
Gr^{F'}(X\ti S^2)\ =\ A_K(t)^2 (1-t)^{2n-4}
\label{last.L1.1}
\end{equation}
where $A_K$ is the Alexander polynomial and $t$ corresponds to the
class $F'$ in $H_2(X\ti S^2)$.
 \label{last.L1}
\end{lemma}

\pf We first show that any class $A\in H_2(X\ti S^2)$ with $RT(A)\neq 0$ is
a multiple of $F'$.  Fix
such a class $A=\alpha\ti \{\mbox{pt}\}+b\{\mbox{pt}\}\ti [S^2]$.   Choose
a sequence of
$({J}_k,{\nu}_k)$-holomorphic maps $f_k$ with $\nu_k\to 0$ and
 $J_k$ converging to a product structure $J\ti j$.  Then the projections
$p\circ f_k$ converge to a $J$-holomorphic bubble tree in $X$ representing
$\alpha$.
But by Theorem \ref{fullGr=Grt} any
$J$-holomorphic curve in $X$ is a multiple of the fiber class $F$.  The
remark following
(\ref{1.defdegzero}) then shows that  $RT(A)= 0$ unless $0=\kappa_{X\ti
S^2}\cdot A$. Using Lemma
\ref{canonical class}, this condition is
\best
0=\kappa_X\cdot \alpha-b\kappa_{ S^2}\cdot S^2=(2g-2+n)F\cdot mF-2b =-2b.
\eest
Thus $A=mF\ti \{\mbox{pt}\}=mF'$ for some $m$.  With this,   Lemma
\ref{GrX*S2} and
 Theorem \ref{fullGr=Grt} give
\best
Gr(X\ti S^2)= Gr^{f}(X\ti S^2)=
 \l[Gr^{F}(X)\r]^2= \l[A(t)(1-t)^{n-2}\r]^2 \;\;\;\Box
\eest

\vskip.2in

 Theorem \ref{7.thm1} now follows easily. When $g>1$, there are
infinitely many monic, integer coefficients, symmetric polynomial  of
degree $2g$  s.t. $A(1)=\pm 1$. For any such
polynomial Burde [B] proved that there is a fibered knot $K$
that has $A$ as its Alexander polynomial.  This gives infinitely many
distinct  Alexander polynomials for
each $g$.   The corresponding spaces $E(n,K)\ti S^2$  then have distinct
Gromov series by Lemma
\ref{last.L1}, so define distinct symplectic structures.


\end{document}